# On the Lottery Problem:
# Tracing Stefan Mandel´s Combinatorial Condensation

Ralph Stömmer *

* Private Researcher

e-mail: ralph.stoemmer@trajectorix.de

**Abstract**

Stefan Mandel is said to have won the lottery 14 times. He never disclosed the recipe he called combinatorial condensation, which enabled him to hit the Romanian lottery jackpot in the early phase of his betting career. Combinatorial condensation is frequently mixed up with another strategy known as buying the pot, which Stefan Mandel was pursuing later on. On occasion, he dropped a few hints on combinatorial condensation. The hints are applied in this work to narrow down and assess his initial recipe. The underlying theory resembles what a weekend mathematician, as he once referred to himself, may have encountered in the 1960s. The cardinality of the (15, 6, 6, 5) and (49, 6, 6, 5) lottery schemes shows that Stefan Mandel probably wasn´t aware of lottery designs. First concepts on such topics had been available at that time, but coherent theories on combinatorial designs took off only in later decades. There is no absolute evidence beyond doubt, but the size of actual covering designs and the associated efforts allow the following hypothesis: Stefan Mandel most likely pioneered in constructing a (15, 6, 6, 5) covering design many years before others published about it, which he applied in the Romanian lottery. Calculations indicate that he took considerable risks that his method might fail, because he limited himself to a subset of 15 numbers out of 49 numbers. The risks explain why he later changed his strategy from combinatorial condensation to buying the pot.

# 1. Introduction

## 1.1. Background on combinatorial minimization

There are many lotteries around the world. A common lottery scheme is (49, 6, 6, t), known as 6/49 lottery, where 6 numbers are randomly drawn by the lottery agency from 49 numbers in an urn without replacement. Beforehand, people buy tickets and choose 6 numbers per ticket. People achieve a t-hit if the 6 numbers drawn intersect the 6 numbers on their ticket in t common elements. The amount of the prize depends on the country, the lottery agency, and the rules of the game. In general, the prize scales with t. For some lotteries the prize starts at $t = 3$, sometimes sufficient to cover the ticket costs, up to $t = 6$ for the jackpot. There are $\binom{49}{6} = 13{,}983{,}816$ combinations to choose 6 from 49 numbers.

Let´s consider four combinations chosen at random, such as {1, 7, 16, 24, 42, 43}, {2, 7, 10, 16, 17, 42}, {5, 7, 16, 28, 33, 42}, {3, 7, 8, 16, 19, 42}. The numbers 7, 16, and 42 appear as common elements. Let p denote the quantity of elements drawn from a total of n elements, and let s, $1 < s \leq p$, denote the quantity of any common elements, then the multitude of combinations with any s common elements MCE(s) in the p/n lottery is given by

$$\mathrm{MCE(s)} = \binom{n-s}{p-s}. \qquad (1)$$

In the 6/49 lottery, with $p = 6$ and $n = 49$, any 3 common elements appear 15,180-fold, 4 common elements appear 990-fold, and 5 common elements appear 44-fold. Any 6 common elements are unique in a 6/49 lottery draw without replacement, they simply appear 1-fold. At an early stage, the question arose if it was possible to reduce the number of combinations required for a t-hit by diminishing the overlap of combinations with common elements.

Documented efforts on related topics date back to the 1850s, when Thomas Kirkman challenged his contemporaries with the “fifteen schoolgirl problem”. Fifteen individuals should be arranged in triples, 5 triples per day on 7 days, in such a way that any two individuals only meet once in a triple [1]. There are $\binom{15}{3} = 455$ combinations to choose 3 from 15 individuals. Equation 1 yields the multitude of combinations that any 2 of them meet more than once, which is $\binom{15-2}{3-2} = 13$. It indicates that 455 combinations could further be reduced by an order of magnitude. First systematic solutions appeared over 100 years after Kirkman formulated the challenge [2, 3], and research on this topic is still ongoing [4].

### 1.2. Status of combinatorial designs during Stefan Mandel´s early activities

When Stefan Mandel started first activities in Romania presumably in the late 1950s and early 1960s, the field of combinatorial designs to tackle lotteries was still a niche. Access to worldwide scientific literature on such topics was probably limited by the Iron Curtain, although a few contemporary experts in the field such as Turan, Sos, or Erdös were researching on the same side [5, 6, 7, 8]. Even if access to their work was easier, there is doubt that Stefan Mandel was building on their achievements in combinatorics, because he never mentioned anyone or any related content in interviews. He referred to the 13th century mathematician Leonardo Fibonacci instead [9]. Lacking scholarly networks where relevant literature circulates, laymen must be able to recognize that topics treated in a different terminology might refer to their problem. If one goes way back in time, subtracts the internet, the telephones at home, the international symposia, and slows down the cross-border letter mail: how should some amateur in Romania have come up with the idea that some limit theorem published in 1963 in Hungary [8], may be helpful for tackling lotteries? Beyond that, research on the lottery problem to reduce the number of combinations was just mastering 2-hits, reaching out to tackle 3-hits. A remark from Hanani et al. [6] sheds light on the situation in 1964:

"If in the lottery problem we want to construct in a similar way a minimal system S which assures a 3-, 4- or 5-hit, we would need a generalization of Turan´s theorem and a generalization of Hanani´s theorem and construction."

Around this time Stefan Mandel devised a number-picking algorithm based on a method he called combinatorial condensation, which would predict five out of six winning numbers [10]. However the method works, Stefan Mandel seems to have found a minimized combinatorial design to tackle lotteries, which he first applied at the Romanian 6/49 lottery, and which he kept a secret ever since. Given the fact that an upper limit of combinations for achieving at least a 3-hit in the 6/49 lottery was still being worked on in the 1970s [11, 12, 13], one may ask: was Stefan Mandel ahead of his time?

### 1.3. Details on Stefan Mandel´s combinatorial condensation

It is not uncommon that amateurs draw level or occasionally are ahead of academics in novel or in niche disciplines. This often proves true for lotteries, where both camps complement and drive each other forward. Mathematicians publish proofs on bounds of designs [14], whereas others, building on computing power and search algorithms, just upload specific constructions of combinatorial designs to online repositories [15, 16, 17].

According to an investigation by reporter Zachary Crockett from The Hustle, Stefan Mandel disclosed a few hints on his combinatorial condensation [9]. In the Romanian 6/49 lottery, he selected 15 numbers, which would have required purchasing $\binom{15}{6} = 5{,}005$ combinations. With 5,005 combinations the odds of winning the jackpot increased to 1 in 2,794. Stefan Mandel claimed his algorithm would reduce these 5,005 combinations to just 569. If the 6 winning numbers fell among the 15 picks, he would be guaranteed to win at least a 2nd prize, that is a 5-hit, and hundreds of smaller prizes, which are 4-hits and 3-hits. Winning the jackpot, a 6-hit, was a 1 in 10 chance. A year later, on the same topic, journalist Juliana Lima from The Sun claimed Stefan Mandel´s combinatorial condensation reduced the before mentioned 5,005 combinations to 669, instead of 569 [18]. His chance to hit the jackpot was 1 in 9, instead of 1 in 10. If Stefan Mandel was wrongly cited, or if he couldn´t come up with the exact details after decades have passed, couldn´t be clarified. For the basic conclusions of this work, the difference doesn´t carry much weight. In the further course, and for the sake of clarity, we will refer to the 569 combinations researched by Zachary Crockett. The 669 combinations from Juliana Lima will be mentioned accordingly.

## 2. Theory

In order to shed light onto the challenges Stefan Mandel had to face, one needs to understand the basic properties of lottery designs and of covering designs. The terminology in this discipline has solidified in the past decades. We follow the notations provided by Bate, van Rees, Colbourn, and Stinson [19, 20, 21]. The definitions for the lottery design and for the covering design are provided first, followed by an example, where both concepts are explained with a simple 3/7 lottery. For reasons of completeness, the example section is supplemented with the logical conclusion and with the probabilistic method.

### 2.1. Definition of lottery design

An (n, k, p, t) lottery scheme consists of an n-set of elements $V = \{1, 2, 3, \ldots, n\}$, a collection of k-element subsets of V (called k-subsets or blocks), and a p-element subset of V (called p-subset). The parameters n, k, p, t must satisfy the relationship $t \leq \{p, k\} \leq n$. The k-subsets, or blocks, are designated as a lottery design LD(n, k, p, t) when any p-subset of V intersects at least one of these blocks in at least t elements. The smallest quantity of blocks of the lottery design LD(n, k, p, t) is denoted as lottery number L(n, k, p, t).

The lottery problem is the fundamental question: what is the smallest quantity of blocks of the lottery design, and what is the lottery design in explicitly written out blocks?

### 2.2. Definition of covering design

A (n, k, t) covering design is a collection of k-element subsets (called k-subsets or blocks) of an n-set of elements V = {1, 2, 3, ..., n} such that any t-element subset (called t-subset) is contained in at least one block. The parameters n, k, t must satisfy the relationship $t \leq k \leq n$. The smallest quantity of blocks of the (n, k, t) covering design is denoted as covering number C(n, k, t).

Covering problems have rarely made headlines, as can be witnessed with lottery problems. However, the challenge to find the smallest quantity of blocks of covering designs, supplemented with explicitly written out blocks, is similar.

Covering designs and lottery designs are related to each other [19]. Adding the trivial case that all k-subsets of an n-set, with gives $N = \binom{n}{k}$ combinations in total, do provide the biggest covering design, L(n, k, p, t) and C(n, k, t) fulfill the inequalities

$$L(n, k, p, t) \leq C(n, k, t) \leq \binom{n}{k} . \tag{2}$$

### 2.3. Example: a 3/7 lottery to illustrate crucial concepts

Let´s assume a (7, 3, 3, 2) lottery scheme, which corresponds to a simple 3/7 lottery, also called the Transylvanian lottery [22]. p = 3 numbers drawn from n = 7 numbers at random without replacement. Beforehand, people buy tickets and choose k = 3 numbers per ticket. A 2-hit is achieved if the p = 3 numbers drawn intersect the k = 3 numbers on at least one of the tickets in t = 2 common elements. All $N = \binom{7}{3} = 35$ combinations are

{1, 2, 3}, {1, 2, 4}, {1, 2, 5}, {1, 2, 6}, {1, 2, 7}, {1, 3, 4}, {1, 3, 5},

{1, 3, 6}, {1, 3, 7}, {1, 4, 5}, {1, 4, 6}, {1, 4, 7}, {1, 5, 6}, {1, 5, 7},

{1, 6, 7}, {2, 3, 4}, {2, 3, 5}, {2, 3, 6}, {2, 3, 7}, {2, 4, 5}, {2, 4, 6},

{2, 4, 7}, {2, 5, 6}, {2, 5, 7}, {2, 6, 7}, {3, 4, 5}, {3, 4, 6}, {3, 4, 7},

{3, 5, 6}, {3, 5, 7}, {3, 6, 7}, {4, 5, 6}, {4, 5, 7}, {4, 6, 7}, {5, 6, 7}.

The question is: how many tickets (alternatively named blocks, or in this case 3-subsets or triples) are required as a minimum to achieve at least a 2-hit?

### 2.3.1. Solution to the 3/7 lottery applying a lottery design

The lottery design LD(7, 3, 3, 2) minimizes the number of tickets for achieving a 2-hit. Any 2 common elements appear $\binom{7-2}{3-2} = 5$-fold, which may be verified by visual inspection of all 35 combinations in section 2.3. The multitude of common elements allows to estimate a first rough reference for L, which is $L \approx 35/5 = 7$. There is no general theory for obtaining efficient lottery designs. For small n, they can be constructed by trial and error.

A lottery design with L(7, 3, 3, 2) = 7 may consist, for example, of the following combinations LD(7, 3, 3, 2) = {{1, 2, 3}, {2, 3, 4}, {3, 4, 5}, {4, 5, 6}, {5, 6, 7}, {6, 7, 1}, {7, 1, 2}}. It is generated, starting from the initial triple {1, 2, 3}, by increasing the elements of each triple {x, y, z} by 1 digit to yield the next triple {x+1, y+1, z+1}. When any element in any triple has reached the maximum number 7, the element in the following triple restarts with the number 1, as done with the last two triples of the lottery design. When the elements of the 7th triple {7, 1, 2} would be used to construct an 8th triple, it would read {1, 2, 3}, which is identical to the first. The 7 blocks fulfill the definition of a lottery design. Any 3-subset of V, provided with the 35 combinations given in section 2.3, intersects at least one of the 7 blocks in the list above in at least 2 elements. Anyone submitting the lottery design would be guaranteed at least a 2-hit. The lottery design is just one of several designs with 7 blocks, which can further be reduced. It is known that LD(7, 3, 3, 2) = 4 [14]. One of several minimized collections of blocks is LD(7, 3, 3, 2) = {{1, 2, 5}, {3, 4, 7}, {2, 5, 6}, {1, 2, 6}}. It fulfills the definition of a lottery design, too. Anyone submitting it to a 3/7 lottery achieves at least a 2-hit.

More sophisticated methods, for higher n, apply computer programs and make use of back-track algorithms, greedy algorithms, simulated annealing, constraint programming and graph theory to reduce the cardinality of blocks [14, 23, 24, 25].

### 2.3.2. Solution to the 3/7 lottery applying a covering design

Covering designs don´t require an additional p-subset to intersect with at least one of the k-subsets in at least t elements. Instead, covering designs are complete in that any t-subset is contained in at least one block. Just as a note, the 2-subset {1, 4} is missing in both lottery designs in the previous section 2.3.1. The 2nd lottery design with LD(7, 3, 3, 2) = 4 is missing even more 2-subsets.

The (7, 3, 2)-covering design with C(7, 3, 2) = 7 consists, for example, of the combinations {1, 2, 6}, {2, 3, 7}, {3, 4, 1}, {4, 5, 2}, {5, 6, 3}, {6, 7, 4}, {7, 1, 5}. It is generated the same way as the first lottery design in section 2.3.1, but starting with the initial triple {1, 2, 6}. The 7 blocks fulfill the definition of a covering design. Any 2-subset is contained in at least one block. Anyone submitting the covering design above to a 3/7 lottery would be guaranteed at least a 2-hit, too. The 7 blocks of the covering design can´t be further reduced without loosing 2-subsets. In contrast to lottery designs, covering designs can be developed more systematically. A simple algorithm to develop a (7, 3, 2) covering design, without guessing around, would be as follows:

i) Record all $N = \binom{7}{3} = 35$ combinations in ascending manner (done in section 2.3).
ii) Start the covering design with a low digit triple, here {1, 2, 3}.
iii) Delete all following triples which contain the 2-subsets {1, 2}, {1, 3}, and {2, 3}.
iv) The first non-deleted triple following {1, 2, 3} is the 2nd candidate for the covering design, here {1, 4, 5}.
v) Delete all following triples which contain the 2-subsets {1, 4}, {1, 5}, and {4, 5}.
vi) The first non-deleted triple following {1, 4, 5} is the 3rd candidate for the covering design, here {1, 6, 7}, and so forth.

The procedure above does lead to a C(7, 3, 2) = 7 covering design with the combinations {1, 2, 3}, {1, 4, 5}, {1, 6, 7}, {2, 4, 6}, {2, 5, 7}, {3, 4, 7}, {3, 5, 6}. The advantage of the approach above is that it is conclusive to the amateur. It can be expanded to find designs of higher order in n, k, and t. The simplicity of the algorithm comes with the disadvantage that considerable computing power is required for higher orders [26].

Research for more sophisticated and faster ways is ongoing, for example by using Fano-planes [25]. Various covering designs, with method of construction, are readily available in the La Jolla Covering Repository [15].

### 2.3.3. Solution to the 3/7 lottery applying logical conclusion

The probability P(t) to achieve a t-hit is given by the fraction $M_t/N$, the classical definition of probability. N denotes the total number of combinations, $N = \binom{n}{p}$, and $M_t$ the number of all possible combinations with a t-hit. It is given by

$$M_t = \binom{p}{t}\binom{n-p}{p-t}. \quad (3)$$

For the 3/7 lottery, the values for $M_t$ , t = 0, 1, 2, 3, are $M_0 = 4$, $M_1 = 18$, $M_2 = 12$, $M_3 = 1$. With respect to the condition

$$\sum_{t=0}^{3} M_t = N \quad (4)$$

one can play it safe for achieving at least a 2-hit by filling out one more block than there are possible combinations $M_0$ with 0-hits and $M_1$ with 1-hits. Figuratively speaking, this method "saturates" possible hits with lower order t = 0, 1. The number of blocks v has to fulfill the condition $v = M_0 + M_1 + 1 = 23$. Rewritten as $v = N - \sum_{t=2}^{3} M_t + 1$, the expression corresponds to a graph theoretic upper bound for the lottery number provided by Gründlingh [14]. The blocks can be filled out at random, but each block has to be unique (blocks have to be pairwise different). The approach guarantees a 2-hit, but filling out 23 blocks, which covers almost 2/3 of all 35 combinations, is not efficient. The method definitely reaches its limits with a (49, 6, 6, 3) lottery scheme, where $v = M_0 + M_1 + M_2 + 1 =$ 13,723,193. Finally, it wouldn´t make sense with (n, 6, 6, 5) lottery schemes for achieving 5-hits, since $v = \sum_{t=0}^{4} M_t + 1 = N - M_5$ .

#### 2.3.4. Solution to the 3/7 lottery applying the probabilistic method

Let us now discuss the strategy in which we pick v blocks at random, but making sure that each block is unique (no doubles, all blocks have to be pairwise different). $M_t$ denotes the number of all possible combinations with a t-hit, as provided in equation 3. $P(m_t, v)$ defines the probability to achieve exactly $m_t$ t-hits, when v blocks are picked at random with no doubles. It is given by the hypergeometric distribution

$$P(m_t, v) = \frac{\binom{M_t}{m_t}\binom{N-M_t}{v-m_t}}{\binom{N}{v}} . \quad (5)$$

$P(m_t = 1, v = 1)$ is the probability to achieve exactly one t-hit with a single block. It is

$$P(m_t = 1, v = 1) = \frac{\binom{M_t}{1}\binom{N-M_t}{0}}{\binom{N}{1}} = \frac{M_t}{N} . \quad (6)$$

The fraction $M_t/N$ on the right hand side in equation 6 is the probability P(t) mentioned in section 2.3.3. Adapted to the 3/7 lottery, the multivariate hypergeometric distribution provides the probability to achieve $m_t$ t-hits, t = 2, 3, with v blocks:

$$P(m_2, m_3, v) = \frac{\binom{M_2}{m_2}\binom{M_3}{m_3}\binom{N-M_2-M_3}{v-m_2-m_3}}{\binom{N}{v}} \ . \tag{7}$$

$P(m_2, m_3, v)$ is required below to calculate the complementary probability for achieving a minimum number of t-hits. At first, we ask for the probability that we achieve no 2-hit and no 3-hit with v blocks. With equation 7 and $N = 35$, $M_2 = 12$, $M_3 = 1$, $m_2 = 0$, $m_3 = 0$, we get

$$P(m_2 = 0, m_3 = 0, v) = \frac{\binom{12}{0}\binom{1}{0}\binom{35-12-1}{v-0-0}}{\binom{35}{v}} = \frac{22!}{35!}\frac{(35-v)!}{(22-v)!} \ . \tag{8}$$

The probability $P(m_2 = 0, m_3 = 0, v)$ is occasionally defined as risk, because one ends up with "nothing". The complement of "nothing" is "at least something", occasionally defined as a level of safety. The complementary probability to achieve at least one 2-hit with v blocks is given by

$$P(m_2 \geq 1, v) = 1 - P(m_2 = 0, m_3 = 0, v) = 1 - \frac{22!}{35!}\frac{(35-v)!}{(22-v)!} \ . \tag{9}$$

To be more specific about the argument in brackets provided in $P(m_2 \geq 1, v)$: achieving at least one 2-hit includes 1 or 2 or 3 … or finally all $M_2 = 12$ combinations with 2-hits or the 3-hit (with $M_3 = 1$ there is only one 3-hit). Equation 9 contains a result previously found by logical conclusion in section 2.3.3. For $v \geq 23$, $P(m_2 \geq 1, v \geq 23) = 100\%$ is a certain event. Table 1 lists a selection of probabilities to achieve at least one 2-hit with the associated number of blocks, calculated with equation 9.

One should internalize the difference between "exactly" one 2-hit, and "at least" one 2-hit. In a 3/7 lottery with $N = 35$ combinations in total, there are at most $M_2 = 12$ combinations with a 2-hit. With equation 5, the probability to achieve exactly one 2-hit ($m_2 = 1$) with $v = 6$ blocks calculates to $P(m_2 = 1, v = 6) = 24.9\%$. In contrast to that, the probability to achieve at least one 2-hit with 6 blocks is calculated with equation 9, which yields $P(m_2 \geq 1, v = 6) = 95.4\%$. The value is considerably higher, it includes other possible – but not explicitly listed - events such as two 2-hits and zero 3-hits up to five 2-hits and one 3-hit.

| (7, 3, 3, 2) lottery scheme | | | | | | | | |
|---|---|---|---|---|---|---|---|---|
| Probability for at least one 2-hit $P(m2 \geq 1, v)$ | 76.5% | 86.0% | 95.4% | 97.5% | 98.6% | 99.3% | 99.9% | 100% |
| Number of blocks v | 3 | 4 | 6 | 7 | 8 | 9 | 12 | 23 |

**Table 1**: Blocks v required to achieve at least one 2-hit for various probabilities $P(m_2 \geq 1, v)$.

With the simple example of the (7, 3, 3, 2) lottery scheme, the advantages and the disadvantages of the probabilistic approach are evident. A construction of lottery designs or of covering designs is not required at all. The only prerequisite is that all blocks, which are picked at random, are unique (no doubles, all blocks have to be pairwise different). The disadvantage is that the quantity of blocks required for acceptable probabilities to achieve at least one 2-hit – for instance $P(m_2 \geq 1, v) > 99\%$ with $v \geq 9$ blocks - is higher compared to lottery designs, where L(7, 3, 3, 2) = 4, or compared to covering designs, where C(7, 3, 2) = 7. Both lottery and covering designs guarantee 100% to achieve at least one 2-hit.

## 3. Results

The statements of Stefan Mandel on combinatorial condensation are related to the probabilistic method, to lottery designs, and to covering designs. The circumstances over 60 years ago will be taken into account.

### 3.1. The limitation to 15 numbers in the 6/49 lottery

Stefan Mandel never disclosed further details why he limited himself right from the start to 15 numbers out of 49 numbers. The limitation basically creates a smaller 6/15 lottery. Ticket costs might have been the primary reason. $N^* = \binom{15}{6} = 5{,}005$ combinations exceeded his budget and that of his companions [9]. Another reason might have been to avoid combinatorial complexity. Lacking computers and even simple electronic calculators, he was restricted to paper and pencil. Whatever he was up to, he first needed to write down all combinations. With time and patience, 5,005 combinations can be listed on 100 pages of squared paper, with ~ 50 combinations per side. All blocks can be generated systematically: each element in the 6-subset has to be changed digit by digit between 1 – 15, avoiding doubles, a dull but straight forward task. The resulting number of all combinations N* provides the maximum quantity of blocks v then applied in the 6/49 lottery. For the 6-hit, equation 5 yields

$$P(m_6 = 1, v = 5{,}005) = \frac{\binom{1}{1}\binom{N-1}{v-1}}{\binom{N}{v}} = \frac{v}{N} = \frac{5{,}005}{13{,}983{,}816} = 0.036\% \,. \qquad (10)$$

The result corresponds with Stefan Mandel´s statement, that 5,005 tickets do increase the odds of winning the jackpot to 1 in 2,794 [9].

In order to quantify Stefan Mandel´s strategy to select 15 picks, we ask: what is the probability, that all 6 winning numbers, or at least 5 of the 6 winning numbers, fall among the n* = 15 picks chosen from n = 49 numbers? The logic behind is equal to asking for the probability that a draw of p = 6 balls from an urn of n = 49 balls, which contains n* = 15 white balls and (n - n*) = 34 black balls, does include t = 5 or t = 6 white balls. The answer is given again by a hypergeometric distribution

$$P(t, n^*) = \frac{\binom{n^*}{t}\binom{n-n^*}{p-t}}{\binom{n}{p}} \ . \tag{11}$$

Equation 11 has the same structure like equation 5, but one should internalize that the perspective has turned to the p = 6 winning numbers. With t = 6, n* = 15, p = 6, and n = 49 equation 11 yields

$$P(t = 6, n^* = 15) = \frac{\binom{15}{6}\binom{34}{0}}{\binom{49}{6}} = \frac{5{,}005}{13{,}983{,}816} = 0.036\% \ . \tag{12}$$

With t = 5, n* = 15, p = 6, and n = 49 equation 11 yields

$$P(t = 5, n^* = 15) = \frac{\binom{15}{5}\binom{34}{1}}{\binom{49}{6}} = \frac{102{,}102}{13{,}983{,}816} = 0.730\% \ . \tag{13}$$

That at least the 5-hit, that is the 5-hit or the 6-hit, falls among the 15 picks is the sum of equations 12 and 13. The sum is provided with fractions containing explicit quantities, which are required later for illustration in section 3.4:

$$P(t \geq 5, n^* = 15) = \frac{102{,}102}{13{,}983{,}816} + \frac{5{,}005}{13{,}983{,}816} = 0.766\% \ . \tag{14}$$

0.766% is a surprisingly low value. Despite the risk associated with limiting the choice to 15 picks, it must have been something else that made Stefan Mandel quite confident. The New Zealand Herald wrote [10]: "… Mr. Mandel boasted he could predict five out of six winning numbers …". He was certain to achieve at least one 5-hit, which is inline with what he told The Hustle [9], where the story continues: "… if the 6 winning numbers fell among the 15 picks…". Stefan Mandel was referring to the 15 picks he selected beforehand, and he referred to his combinatorial condensation, which he claimed would reduce the number of blocks from 5,005 down to 569 [9]. How did he make sure to have secured at least one 5-hit?

### 3.2. The LD(15, 6, 6, 5) and LD(49, 6, 6, 5) lottery designs

Up to this point, the probability theory helped to quantify the risk Stefan Mandel took upon himself by limiting his game to 15 picks out of 49 numbers, and the value for the 6-hit fits in with the statements to the press. Why did he spend considerable time and efforts to pore over combinations?

It is unlikely that he developed a lottery design for a (15, 6, 6, 5) lottery scheme, not to mention for a (49, 6, 6, 5) lottery scheme. As noted in section 1.2, scientists in the early 1960s started to tackle 3-hits [6]. Higher order lottery designs for 5-hits require computing power, which was unavailable at that time. Recent activities reveal the challenges to construct LD(15, 6, 6, 5) and LD(49, 6, 6, 5) lottery designs. Before explicit designs came up, the lower and the upper bounds for the lottery numbers L(15, 6, 6, 5) and L(49, 6, 6, 5) were determined to at least narrow down the challenge. In this context, the work of Gründlingh is significant, because it provides a survey of achievements up to 2004 [14]. It seems that addressing 5-hits didn´t come up long before that date, because in 1996 Colbourn and Dinitz edited a book on the state of the art of combinatorial designs, where the lottery schemes discussed were still limited to 3-hits [20]. Obviously, better performing computing hardware and software developments considerably enhanced activities. From that, one can conclude that information on higher order lottery schemes was not available when Stefan Mandel developed his combinatorial condensation. Explicit lists with blocks appeared much later, made accessible in online repositories. Table 2 provides a brief overview for the two lottery schemes in question. Where possible, the earlier dates of first uploads are provided, too. Stefan Mandel´s combinatorial condensation reduced 5,005 combinations to 569 blocks, which is not sufficient for a lottery design based on 49 numbers. Lacking computing hardware, he was unable to come up with designs with block numbers going into the thousands. This challenge is currently being addressed, as the dates in table 2 demonstrate.

It can be excluded that designs of such scale have been around unpublished for long for secretly testing them out. Bulk purchasing of lottery tickets as well as setting up syndicates for submitting thousands of tickets do violate rules of many lotteries, and a considerable amount of money and logistics for submitting tickets would be required. In this regard, restrictions to lottery regulations were implemented by many agencies shortly after Stefan Mandel and his syndicate scooped several jackpots.

Things do seem less complex for the LD(15, 6, 6, 5) lottery design. Given the lower bound for the lottery number, a lottery design was probably on the horizon. A lottery design with 142 blocks is online with a time stamp 2018, a first upload dates back to 2006. However, we can exclude that Stefan Mandel created a lottery design, because his algorithm did arrive at 569 blocks, which is a factor of 4 above actual lottery design proposals.

| | Lottery number | | Specific lottery design | |
|---|---|---|---|---|
| | lower bound | published | blocksize | published |
| (15, 6, 6, 5) lottery scheme | 100 [14] | 2004 | 143 [17]<br>142 [16, 17] | 2006<br>2018 |
| (49, 6, 6, 5) lottery scheme | 62,151 [14] | 2004 | 151,012 [17]<br>144,819 [17]<br>142,361 [16] | 2012<br>2023<br>2024 |

**Table 2**: Basic data on lottery numbers and lottery designs for 2 lottery schemes.

### 3.3. The (15, 6, 5) covering design

It was mentioned in section 2.2 that $C(n, k, t) \geq L(n, k, p, t)$, so we dispense with the discussion if Stefan Mandel might have created a (49, 6, 5) covering design decades before sufficient computing power was available. Table 3 provides basic data on the (15, 6, 5) covering design, framed with nearby designs. A first upload dates back to 1996, with the number of blocks reduced from 620 down to 578 in 2007. Referring to LD(15, 6, 6, 5) in table 2, the covering design was available earlier. Without overestimating the time stamps of the uploads, the decade in between covering design and lottery design indicates that covering designs are constructed easier in a straight forward manner, as mentioned in section 2.3.2.

| | Covering number | | Specific covering design | |
|---|---|---|---|---|
| | lower bound | published | blocksize | published |
| (14, 6, 5) covering design | 348 [15] | 1996 | 385 [15]<br>371 [15] | 1996<br>1997 |
| (15, 6, 5) covering design | 548 [15] | 1996 | 620 [15]<br>578 [15] | 1996<br>2007 |
| (16, 6, 5) covering design | 731 [15] | 1996 | 840 [15]<br>808 [15] | 1996<br>1996 |

**Table 3**: Basic data on covering numbers and covering designs for relevant coverings.

Nobody knows if Stefan Mandel was considering something like equation 1 to anticipate that he could reduce 5,005 combinations by roughly one order of magnitude, since any 5 common elements do appear 10-fold in a 6/15 lottery. However he got there, his 569 blocks correspond with the values of the (15, 6, 5) covering design in table 3. With respect to the conditions 60 years ago, the hypothesis is as follows: Equipped with paper and pencil, he first might have recorded all $N^* = \binom{15}{6} = 5{,}005$ combinations in ascending manner. Then he moved forward with an algorithm similar to that one outlined in section 2.3.2. He might have started with a low digit 6-subset, such as {1, 2, 3, 4, 5, 6}, the 1st candidate for the covering design. Then he deleted all following combinations which contained the 5-subsets {1, 2, 3, 4, 5}, {1, 2, 3, 4, 6}, {1, 2, 3, 5, 6}, {1, 2, 4, 5, 6}, {1, 3, 4, 5, 6}, {2, 3, 4, 5, 6}. The first non-deleted combination then was the 2nd candidate for the covering design, which in turn had to be split into 5-subsets, and so forth. Deleted combinations were out of the game, so for each search-and-delete loop, the quantity of combinations to be checked for specific 5-subsets decreased.

We can only roughly estimate the complexity of the task, when Stefan Mandel checked all combinations downstream if they contained any of the six 5-subsets from the 1st to the final 569th candidate of the covering design. At the beginning of the procedure one doesn´t know the final result for the covering number. In retrospect, 569 blocks with six 5-subsets each yield 569 x 6 = 3,414 5-subsets in total. The 3,414 5-subsets do comprise the loops the search space, which reduces after each loop, is subjected to. We don´t know if and how often he was back checking deleted blocks upstream to be on the safe side. As a rough rule of thumb, let´s assume that the search space decreases linearly from the first candidate of the covering design on top of the list of combinations to the final candidate at the end. With great margin, the number of checks goes into several millions. With undisturbed focus, on average 4 s are needed for checking and deleting combinations, verified in a self-test by the author on a few pages filled with combinations. Scaling up the self-test, 1 million checks would result in ~ 1,111 h. Given 4 h per day spare time, the work would result in 278 days work, without leisure. With leisure, 1 million checks would roughly amount to a year. Basically, the whole enterprise was a question of accuracy, patience, and time. Time is what Stefan Mandel probably had, because the New Zealand Herald cites [10]: "… after years of painstaking research, he devised a number-picking algorithm, which drew on a method he called "combinatorial condensation". The Sun refers to the same topic with [18]: "After years of planning and mathematical research".

We do not know how often he tested his combinatorial condensation on fewer or on more picks. Table 3 provides examples for the (14, 6, 5) and the (16, 6, 5) covering design Stefan Mandel may have encountered on his quest. Just 14 picks yield $\binom{14}{6} = 3{,}003$ combinations, and 16 picks yield $\binom{16}{6} = 8{,}008$ combinations. Actual (14, 6, 5) covering designs consist of 371 blocks, and (16, 6, 5) covering designs of 808 blocks, as listed in table 3. 371 blocks with six 5-subsets each yield 2,226 5-subsets, which do comprise the loops the search space of initially 3,003 combinations is subjected to. 808 blocks with six 5-subsets each yield 4,848 5-subsets, which do comprise the loops the search space of initially 8,008 combinations is subjected to. From that, the time required to construct a (14, 6, 5) covering design should be roughly a factor of ~0.5 lower than the time for a (15, 6, 5) covering design, and the time required for a (16, 6, 5) covering design should be roughly a factor of ~2 higher than the time for a (15, 6, 5) covering design. Experiments like this add considerable time budgets, so Stefan Mandel must have carefully considered quantities. In later decades, when computing time was a precious asset, articles appeared on the efforts of search algorithms related to coverings [27]. Li [26] remarks that the search space tends to grow exponentially with n and k. Gordon et al. [28] state that the algorithm would asymptotically take time and space $O(n^k)$. With the latter estimate, a (14, 6, 5) covering design would require a factor of ~0.7 less time, and the (16, 6, 5) covering design would require a factor of ~1.5 more time compared to the (15, 6, 5) covering design. Whatever applies: for Stefan Mandel, time was most likely the limiting factor.

The appendix contains an example of a (15, 6, 5) covering design from the La Jolla Covering Repository, provided by Barnabas Toth [15]. Although the covering number doesn´t match exactly, which may be attributed to the fact that Stefan Mandel couldn´t come up with the exact number after decades have passed, his design should have looked similar.

### 3.4. Summary of probabilities in the 6/49 lottery

Tables 4a and 4b summarize the probabilities to achieve at least one 5-hit or the 6-hit under certain conditions, for two relevant lottery schemes. Table 4a provides the probabilities for Stefan Mandel´s strategy to focus on 15 picks from 49 numbers, whereas table 4b provides the probabilities with blocks chosen completely at random from all 49 numbers.

If all 5,005 blocks resulting from the reduced set of 15 picks were filled out – which Stefan Mandel and his companions couldn´t afford – the probabilities for the smaller 6/15 lottery would be 100% each, as provided in table 4a. As a result of this, the probabilities for the 6/49 lottery are just determined by the restriction to 15 picks, as calculated with equations 12 to 14. Things look slightly different with a (15, 6, 5) covering design. Calculations involving the covering design are carried through with a block size of 569 as provided in reference [9] (the conclusions are similar if the block size of 669 was applied from reference [18]). For the 5-hit, the basic principle of covering designs is confirmed: when at least 5 winning numbers drawn from 49 numbers fall among the 15 picks – for this event the probability is 0.766% - the covering design makes sure that these 5 numbers are covered by at least one block. The probability does remain 0.766% without change. Under these conditions, the covering design does allow one to move into a 6/49 lottery with a nearly 10-fold reduction of the number of blocks, what Stefan Mandel presumably wanted to achieve.

For the 6-hit, the situation is different. One of the 569 blocks might contain the 6-hit. Equation 5 yields for the 6/15 lottery

$$P(m_6 = 1, v = 569) = \frac{569}{5{,}005} = 11.369\% \,. \tag{15}$$

The result corresponds to a note that "… if the 6 winning numbers fell among his 15 picks … he´d have a 1 in 10 chance of winning the grand prize." [9]. The same logic applies to the 6/49 lottery, where equation 5 yields

$$P(m_6 = 1, v = 569) = \frac{569}{13{,}983{,}816} = 0.004\% \,. \tag{16}$$

With respect to the low probability of the 6-hit in the 6/49 lottery, the (15, 6, 5) covering design creates a disadvantage. If Stefan Mandel would have resorted to the brute force approach with 5,005 combinations, the probability to achieve the 6-hit would increase almost 10-fold to 0.036% - but which is still a small value. It looks like a matter of balancing a negligible probability versus a small probability, and Stefan Mandel most likely didn´t care.

Without belittling Stefan Mandel´s efforts with combinatorial condensation, calculations show that he probably might have been slightly better off without. Table 4b provides the probabilities in the 6/49 lottery for choosing blocks completely at random ("completely at random" includes the simplification, that blocks, which happen to be double, are allowed). Achieving at least one 5-hit and the 6-hit with v blocks picked completely at random is given by the probabilities

$$P(m_5 \geq 1, v) = 1 - [1 - P(5) - P(6)]^v \,, \tag{17}$$

$$P(m_6 = 1, v) = 1 - [1 - P(6)]^v \,. \tag{18}$$

With $P(5) = M_5/N$, $P(6) = M_6/N$, $M_5 = 258$, $M_6 = 1$, $N = \binom{49}{6} = 13{,}983{,}816$, equations 17 and 18 are rewritten as

$$P(m_5 \geq 1, v) = 1 - \left[1 - \frac{258+1}{13{,}983{,}816}\right]^v , \tag{19}$$

$$P(m_6 = 1, v) = 1 - \left[1 - \frac{1}{13{,}983{,}816}\right]^v . \tag{20}$$

The corresponding probabilities for v = 5,005 and v = 569 blocks are listed in table 4b. Compared to the values in table 4a, the probabilities for achieving a 6-hit don´t change. Whether one takes all 5,005 blocks resulting from a reduced set of 15 picks from 49 numbers, or one chooses 5,005 blocks completely at random from 49 numbers, the probability for the 6-hit is 0.036%. Whether one constructs a covering design with C(15, 6, 5) = 569 from a reduced set of 15 picks from 49 numbers, or one chooses 569 blocks completely at random from 49 numbers, the probability for the 6-hit is 0.004%. In equations 10, 12 and 20 the probability to achieve a 6-hit scales with the block size v, whether filled out at random or with a strategy (for $N \rightarrow \infty$, equation 20 can be approximated as $P(m_6 = 1, v) \approx v/N$).

The results differ for achieving at least one 5-hit. When 5,005 blocks are chosen at random from all 49 numbers, the probability calculates to 8.853%. If just 569 blocks are chosen at random, the probability is 1.048%. Both probabilities for the random choice in table 4b are higher compared to the corresponding values table 4a, the first one by at least one order of magnitude. A catchy explanation for the difference is challenging. It can be made plausible by comparing equations 14 and 19. Both equations are inherently clear when each one is taken by itself. But to come to an insight through comparison, equation 19 is better presented as a sum of fractions. For $N > (M_5+M_6)$ and v = 5,005, it can be roughly approximated with

$$P(m_5 \geq 1, v = 5{,}005) =$$

$$= 1 - \left[1 - \frac{258+1}{13{,}983{,}816}\right]^{5{,}005} \approx 1 - \left[1 - 5{,}005\,\frac{258+1}{13{,}983{,}816}\right] = \frac{5{,}005 \cdot 258}{13{,}983{,}816} + \frac{5{,}005}{13{,}983{,}816} . \tag{21}$$

The 2nd fraction in equation 14 corresponds to the 2nd fraction in approximation 21. The 1st fraction in equation 14 contains the 102,102 favorable events that 5 winning numbers fall among the 15 picks. Although approximation 21 provides a probability, its 1st fraction resembles the structure of an expected value $\sim v \cdot P(5)$ for 5-hits, when all $M_5 = 258$ combinations for 5-hits in a 6/49 lottery are available. Figuratively speaking, it is better to randomly spread the blocks over the whole space of combinations than to cover just a limited space. The same considerations apply to the smaller quantity of 569 blocks.

| | 15 picks from 49 numbers | | | |
|---|---|---|---|---|
| | All blocks<br>v = 5,005 | Covering design<br>C (15, 6, 5) = 569 | All blocks<br>v = 5,005 | Covering design<br>C(15, 6, 5) = 569 |
| | Probability for at least one 5-hit | | Probability for the 6-hit | |
| (15, 6, 6, 5)<br>lottery scheme | 100% | 100% | 100% | 11.369% |
| Winning numbers<br>among 15 picks | 0.766% | | 0.036% | |
| (49, 6, 6, 5)<br>lottery scheme | 0.766% | 0.766% | 0.036% | 0.004% |

**Table 4a**: Probabilities to achieve at least one 5-hit or the 6-hit for 2 lottery schemes in case all 5,005 combinations of 15 picks are chosen, or if a covering design with 569 blocks is chosen.

| | All 49 numbers | | | |
|---|---|---|---|---|
| | Random blocks<br>v = 5,005 | Random blocks<br>v = 569 | Random blocks<br>v = 5,005 | Random blocks<br>v = 569 |
| | Probability for at least one 5-hit | | Probability for the 6-hit | |
| (49, 6, 6, 5)<br>lottery scheme | 8.853% | 1.048% | 0.036% | 0.004% |

**Table 4b**: Probabilities to achieve at least one 5-hit or the 6-hit for the 6/49 lottery with blocks picked completely at random.

A more striking example which reveals misleading imagination would be to take 25 picks from 49 numbers, which amounts to roughly 50% of all numbers. One could naively assume $P(t \geq 5, n^* = 25) \approx 50\%$, that a fair distribution of luck would scale with the quantity of numbers covered. With $n^* = 25$ inserted into equations 12 and 13, the result emphasizes once again that limiting lottery numbers is generally a bad idea:

$$P(t \geq 5, n^* = 25) = \frac{\binom{25}{5}\binom{24}{1}}{\binom{49}{6}} + \frac{\binom{25}{6}\binom{24}{0}}{\binom{49}{6}} = \frac{1{,}275{,}120}{13{,}983{,}816} + \frac{177{,}100}{13{,}983{,}816} = 10{,}4\,\% \,. \tag{22}$$

However, for Stefan Mandel the restriction to 15 picks went off lightly. As the values in table 4a and 4b show, the covering design assures a probability of 0.766% to achieve at least one 5-hit, whereas 569 blocks chosen at random would assure a probability of 1.048%, which isn´t significantly better on this small scale. Even if he had followed the random path and saved a lot of time, his risk to fail would still have been considerably high. It is a pity that the restriction to 15 picks limits the positive effect of the covering design.

In retrospect, we can only speculate as to why Stefan Mandel put so much effort into combinations. Combinatorics and probability theory are two sides of the same medal. Those who lose themselves on one side may not see the other. A possible explanation would be that end of the 1950s, when Stefan Mandel started his activities in Romania, he focused on combinations because he simply was lacking a calculator. This is hard to imagine from today's perspective, but it is certainly conceivable. Repeated computations of the type as equations 11 or 17, which are crucial to decide which way to advance, are not easily accomplished with paper and pencil, even if a slide rule was available. This includes preparatory computations with varying parameters in several binomial coefficients. Perhaps the tedious task of systematically writing down thousands of combinations and searching them through for common subsets simply did appear more manageable.

### 3.5. Did Stefan Mandel apply a greedy algorithm years before it was circulated?

There is reason to assume that Stefan Mandel was working along a procedure now circulated as "greedy algorithm". First of all, the easiest way to record all $N^* = \binom{15}{6} = 5{,}005$ combinations in an ascending manner would be as accomplished for the smaller triples in section 2.3. One would start with the 6-subset {1, 2, 3, 4, 5, 6} and continue by increasing the 6th digit by +1 until the final number 15 is reached. Then increasing the 5th digit by +1, and continue with the 6th digit as done before, and so forth, until the 1st digit is reached. Experts would call the resulting sequence of 6-subsets a lexicographic order. An algorithm for generating lexicographic orders of k-subsets of an n-set is provided in [29]. Counting up digits can be accomplished with paper and pencil, and one doesn´t lose track. It would be the natural thing to do, even without in-depth mathematical knowledge.

As demonstrated in section 2.3.2, common subsets can be identified and deleted downstream. The procedure follows what is called a "greedy algorithm". We quote the crucial passage from Gordon et al. [28]:

i) Arrange the k-subsets of a n-set in a list.
ii) Choose from the list the k-subset that contains the maximum number of t-sets that are still uncovered. In case of ties, choose the k-subset occurring earliest in the list.
iii) Repeat step ii) until all t-sets are covered

The greedy algorithm is self-explanatory, the task is simple but quite monotonous, which makes it a question of accuracy, patience and time. We didn´t investigate in detail for how long lexicographic orders and greedy algorithms have been around in scientific literature. Their origins are hard to trace back, but it is reasonable to assume that such topics received wider dissemination in the 1970s, when researchers gained frequent access to mainframe computers [27, 30, 31, 32]. There is no indication that Stefan Mandel was drawing from existing sources for greedy algorithms, if any were available to him at that time. However, it would be plausible to assume that he worked through all 5,005 combinations with search-and-delete loops related to a greedy process. Lacking computers, alternatives are hard to imagine which allow to reduce overlaps of common elements as systematically. The associated efforts would fit to Stefan Mandel´s time-budget, and the resulting order of magnitude for the covering number makes sense.

## 4. Conclusion

This work traces Stefan Mandel´s combinatorial condensation, which he applied in the 1960s to hit the jackpot in the Romanian 6/49 lottery. He most likely constructed a first (15, 6, 5) covering design with a block size of 569, according to another source with a block size of 669, years before covering numbers of higher order were calculated to narrow down the challenge, and many years before explicit blocks for higher order covering designs were published.

Section 3.1 shows that he took a considerable risk by limiting the numbers to 15 picks to save ticket costs. Given the circumstances 60 years ago, lacking any computing power, the measure probably also served to reduce the combinatorial complexity to manageable levels. Section 3.4. shows that a (15, 6, 5) covering design allows one to move into a 6/49 lottery with a nearly 10-fold reduction of the number of blocks. However, further calculations reveal that Stefan Mandel could just as well have chosen 569 blocks at random. The 15 picks from 49 numbers determine the probability to achieve at least one 5-hit and limit the positive effect of the covering design. Whatever path he would have taken, it turns out that at the end he was lucky, a circumstance which doesn´t diminish his achievements.

For whatever reason Stefan Mandel never revealed the details of his algorithm, it is unlikely that he just made up the whole story. None of his companions in Romania and elsewhere, with whom he shared the costs for lottery tickets and efforts for submission, did object his story, which one would expect given the worldwide headlines his campaigns made. The hints he gave on several occasions fit together, the calculations confirm the figures, the odd number of 569 blocks (according to another source 669 blocks) he brought up fits well to actual boundaries, and the efforts involved to create a (15, 6, 5) covering design with paper and pencil go with his remarks that the project on combinatorial condensation took years.

Further developments indicate that Stefan Mandel was aware of the risks his method implied, and he switched strategy to buying the pot. Buying the pot requires bulk-purchasing of lottery tickets to cover all combinations, which reduces the risk to zero. The money required to buy all tickets would be raised by a syndicate [33]. In the 1980s affordable personal computers came to the market. Efforts became manageable to program software for printing out all 13,983,816 combinations for the 6/49 lottery, or similar orders of magnitude for other lottery schemes. The main criterion was that the jackpot climbed to levels higher than all costs related to the venture, which includes ticket costs for all combinations, transaction efforts, and costs associated with running an enterprise. This is how Stefan Mandel continued in Australia and later in the USA, at times when bulk-purchasing and submission of tickets by syndicates were largely legal [34].

One wonders why he didn´t expand his combinatorial condensation to higher order lottery schemes when computers were available. Maybe he tried. But given the computing power in the 1980s and early 1990s, it seemed less complex to simply print out all combinations of a lottery, instead of devising sophisticated algorithms. It is an irony of time that an amateur, "an accountant without too much education" as he called himself [9], discarded combinatorial condensation at a moment when research on that topic unfolded momentum.

## Acknowledgment

The author would like to express deep gratitude to the anonymous reviewer for scrutinizing the manuscript, which has gained a lot from his remarks.

# Appendix

Example of a (15, 6, 5) covering design provided by Barnabas Toth [15], C(15, 6, 5) = 578.

| | | | | | | | |
|---|---|---|---|---|---|---|---|
| 1 2 3 4 5 14 | 1 2 7 11 12 15 | 1 4 5 7 9 10 | 1 7 8 9 13 15 | 2 4 5 6 11 13 | 2 6 9 12 13 15 | 3 5 9 10 11 14 | 4 7 8 9 10 15 |
| 1 2 3 4 6 11 | 1 2 8 9 11 13 | 1 4 5 7 12 15 | 1 7 8 10 12 13 | 2 4 5 6 12 14 | 2 6 10 11 12 15 | 3 5 10 11 12 13 | 4 7 8 9 14 15 |
| 1 2 3 4 7 9 | 1 2 8 9 12 13 | 1 4 5 8 9 13 | 1 7 9 10 12 15 | 2 4 5 7 8 15 | 2 6 11 13 14 15 | 3 5 10 12 14 15 | 4 7 8 10 11 15 |
| 1 2 3 4 8 12 | 1 2 8 10 12 15 | 1 4 5 8 10 12 | 1 7 9 10 14 15 | 2 4 5 7 10 13 | 2 7 8 9 10 12 | 3 6 7 8 10 11 | 4 7 8 10 13 14 |
| 1 2 3 4 10 15 | 1 2 8 10 13 14 | 1 4 5 8 11 15 | 1 7 9 11 12 13 | 2 4 5 7 13 14 | 2 7 8 10 11 13 | 3 6 7 8 13 14 | 4 7 8 11 12 13 |
| 1 2 3 4 11 13 | 1 2 8 12 14 15 | 1 4 5 9 10 15 | 1 7 10 11 12 13 | 2 4 5 8 9 12 | 2 7 8 11 13 14 | 3 6 7 9 14 15 | 4 7 8 12 13 15 |
| 1 2 3 5 6 10 | 1 2 9 10 12 14 | 1 4 5 11 12 13 | 1 7 12 13 14 15 | 2 4 5 8 11 14 | 2 7 9 10 13 14 | 3 6 7 10 12 13 | 4 7 9 10 12 13 |
| 1 2 3 5 7 8 | 1 2 9 13 14 15 | 1 4 5 13 14 15 | 1 8 9 10 11 15 | 2 4 5 9 11 12 | 2 7 9 11 14 15 | 3 6 7 10 12 14 | 4 7 9 11 12 15 |
| 1 2 3 5 8 15 | 1 2 11 12 13 14 | 1 4 6 7 8 10 | 1 8 9 10 13 14 | 2 4 5 9 13 15 | 2 7 10 12 13 14 | 3 6 7 10 12 15 | 4 7 9 11 13 14 |
| 1 2 3 5 9 11 | 1 3 4 5 6 8 | 1 4 6 7 9 15 | 1 8 9 11 12 14 | 2 4 5 10 12 15 | 2 7 10 13 14 15 | 3 6 7 11 13 15 | 4 8 9 10 11 14 |
| 1 2 3 5 12 13 | 1 3 4 5 7 12 | 1 4 6 7 11 12 | 1 8 11 13 14 15 | 2 4 5 11 14 15 | 2 8 9 10 13 15 | 3 6 8 9 10 12 | 4 8 9 10 12 14 |
| 1 2 3 6 7 15 | 1 3 4 5 7 15 | 1 4 6 8 9 11 | 1 9 11 12 13 15 | 2 4 6 7 9 14 | 2 8 9 11 14 15 | 3 6 8 9 11 14 | 4 8 9 11 13 15 |
| 1 2 3 6 8 14 | 1 3 4 5 9 11 | 1 4 6 8 12 14 | 1 10 12 13 14 15 | 2 4 6 7 11 15 | 2 8 11 12 13 15 | 3 6 8 9 13 15 | 4 9 10 13 14 15 |
| 1 2 3 6 9 15 | 1 3 4 5 10 13 | 1 4 6 8 13 15 | 2 3 4 5 6 15 | 2 4 6 8 11 13 | 2 9 10 11 12 13 | 3 6 8 11 12 13 | 4 10 11 12 13 15 |
| 1 2 3 6 12 15 | 1 3 4 5 12 15 | 1 4 6 9 10 14 | 2 3 4 5 7 9 | 2 4 6 8 12 15 | 2 9 10 11 14 15 | 3 6 8 12 14 15 | 4 11 12 13 14 15 |
| 1 2 3 6 13 15 | 1 3 4 6 7 14 | 1 4 6 9 12 13 | 2 3 4 5 8 10 | 2 4 6 8 13 14 | 2 9 11 12 14 15 | 3 6 9 10 11 15 | 5 6 7 8 10 12 |
| 1 2 3 7 10 13 | 1 3 4 6 9 12 | 1 4 6 10 12 15 | 2 3 4 5 10 11 | 2 4 6 9 11 13 | 3 4 5 6 7 10 | 3 6 9 10 13 14 | 5 6 7 8 12 14 |
| 1 2 3 7 11 12 | 1 3 4 6 10 13 | 1 4 6 11 13 15 | 2 3 4 5 12 13 | 2 4 6 9 14 15 | 3 4 5 6 8 11 | 3 6 9 12 13 14 | 5 6 7 9 10 15 |
| 1 2 3 7 13 14 | 1 3 4 6 14 15 | 1 4 7 8 9 12 | 2 3 4 6 7 13 | 2 4 6 10 11 14 | 3 4 5 6 8 12 | 3 6 11 12 14 15 | 5 6 7 9 13 14 |
| 1 2 3 8 9 13 | 1 3 4 7 8 13 | 1 4 7 8 11 15 | 2 3 4 6 8 9 | 2 4 6 10 13 15 | 3 4 5 6 8 13 | 3 7 8 9 10 13 | 5 6 7 10 11 13 |
| 1 2 3 8 10 11 | 1 3 4 7 10 11 | 1 4 7 9 11 14 | 2 3 4 6 10 12 | 2 4 6 11 12 13 | 3 4 5 6 9 14 | 3 7 8 9 12 15 | 5 6 7 11 12 15 |
| 1 2 3 9 10 14 | 1 3 4 8 9 10 | 1 4 7 9 13 14 | 2 3 4 6 12 14 | 2 4 7 8 9 13 | 3 4 5 7 8 14 | 3 7 8 11 13 15 | 5 6 8 9 10 11 |
| 1 2 3 9 12 14 | 1 3 4 8 9 14 | 1 4 7 10 12 14 | 2 3 4 7 8 11 | 2 4 7 8 12 14 | 3 4 5 7 11 13 | 3 7 9 10 11 12 | 5 6 8 9 12 13 |
| 1 2 3 10 11 15 | 1 3 4 8 11 15 | 1 4 7 10 13 15 | 2 3 4 7 9 12 | 2 4 7 9 10 11 | 3 4 5 7 12 15 | 3 7 9 11 12 14 | 5 6 8 9 14 15 |
| 1 2 3 10 12 14 | 1 3 4 9 11 15 | 1 4 7 11 13 14 | 2 3 4 7 9 15 | 2 4 7 10 12 15 | 3 4 5 8 9 15 | 3 7 9 12 13 15 | 5 6 8 11 12 15 |
| 1 2 3 11 14 15 | 1 3 4 9 13 15 | 1 4 8 10 11 13 | 2 3 4 7 10 14 | 2 4 7 11 12 14 | 3 4 5 9 10 12 | 3 7 9 13 14 15 | 5 6 8 11 13 14 |
| 1 2 4 5 6 15 | 1 3 4 10 12 13 | 1 4 8 10 14 15 | 2 3 4 8 13 15 | 2 4 7 11 13 15 | 3 4 5 9 13 14 | 3 7 10 11 13 14 | 5 6 9 11 12 14 |
| 1 2 4 5 7 11 | 1 3 4 10 13 14 | 1 4 8 11 12 15 | 2 3 4 8 14 15 | 2 4 8 9 10 14 | 3 4 5 10 14 15 | 3 7 11 12 13 15 | 5 6 10 11 14 15 |
| 1 2 4 5 8 13 | 1 3 4 11 12 14 | 1 4 8 12 13 14 | 2 3 4 9 10 13 | 2 4 8 9 11 12 | 3 4 5 11 12 14 | 3 8 9 10 14 15 | 5 6 10 12 13 14 |
| 1 2 4 5 9 14 | 1 3 5 6 7 13 | 1 4 9 10 11 12 | 2 3 4 9 11 14 | 2 4 8 10 11 15 | 3 4 5 11 13 15 | 3 8 9 11 12 15 | 5 6 11 12 13 15 |
| 1 2 4 5 10 14 | 1 3 5 6 9 10 | 1 4 9 11 13 14 | 2 3 4 11 12 15 | 2 4 8 10 12 13 | 3 4 6 7 8 15 | 3 8 9 11 13 14 | 5 7 8 9 10 13 |
| 1 2 4 5 12 14 | 1 3 5 6 11 15 | 1 4 9 12 14 15 | 2 3 4 13 14 15 | 2 4 9 10 12 15 | 3 4 6 7 9 13 | 3 8 10 12 13 14 | 5 7 8 9 11 14 |
| 1 2 4 6 7 8 | 1 3 5 6 12 14 | 1 4 10 11 14 15 | 2 3 5 6 7 11 | 2 4 9 12 13 14 | 3 4 6 7 11 12 | 3 9 10 12 13 15 | 5 7 8 13 14 15 |
| 1 2 4 6 7 10 | 1 3 5 7 9 13 | 1 5 6 7 8 11 | 2 3 5 6 8 9 | 2 4 10 11 13 14 | 3 4 6 8 10 14 | 3 10 11 13 14 15 | 5 7 9 10 12 14 |
| 1 2 4 6 9 12 | 1 3 5 7 9 14 | 1 5 6 7 10 14 | 2 3 5 6 11 14 | 2 4 10 12 14 15 | 3 4 6 9 10 11 | 4 5 6 7 8 9 | 5 7 9 11 12 13 |
| 1 2 4 6 13 14 | 1 3 5 7 10 11 | 1 5 6 7 12 13 | 2 3 5 6 12 13 | 2 5 6 7 8 13 | 3 4 6 9 10 15 | 4 5 6 7 9 11 | 5 7 9 12 14 15 |
| 1 2 4 7 8 10 | 1 3 5 8 9 11 | 1 5 6 7 13 15 | 2 3 5 7 10 12 | 2 5 6 7 8 14 | 3 4 6 10 11 15 | 4 5 6 7 14 15 | 5 7 10 12 13 15 |
| 1 2 4 7 12 13 | 1 3 5 8 10 12 | 1 5 6 8 9 15 | 2 3 5 7 11 14 | 2 5 6 7 10 15 | 3 4 6 11 13 14 | 4 5 6 8 10 15 | 5 7 11 13 14 15 |
| 1 2 4 7 14 15 | 1 3 5 8 13 14 | 1 5 6 8 10 13 | 2 3 5 7 13 15 | 2 5 6 8 10 11 | 3 4 6 12 13 15 | 4 5 6 9 13 15 | 5 8 9 10 12 15 |
| 1 2 4 8 9 15 | 1 3 5 9 11 12 | 1 5 6 9 11 13 | 2 3 5 8 11 13 | 2 5 6 8 13 15 | 3 4 7 8 9 11 | 4 5 6 10 12 13 | 5 8 9 11 13 15 |
| 1 2 4 8 11 14 | 1 3 5 9 13 15 | 1 5 6 9 14 15 | 2 3 5 8 12 14 | 2 5 6 9 10 12 | 3 4 7 8 10 12 | 4 5 6 10 13 14 | 5 8 10 11 13 14 |
| 1 2 4 9 10 13 | 1 3 5 10 14 15 | 1 5 6 10 12 15 | 2 3 5 9 10 15 | 2 5 6 9 10 13 | 3 4 7 9 10 14 | 4 5 6 11 12 15 | 5 8 11 12 13 14 |
| 1 2 4 9 11 15 | 1 3 5 11 13 14 | 1 5 7 8 9 12 | 2 3 5 9 12 13 | 2 5 6 9 10 14 | 3 4 7 10 13 15 | 4 5 7 8 10 11 | 5 9 10 11 12 15 |
| 1 2 4 10 11 12 | 1 3 6 7 8 9 | 1 5 7 8 10 15 | 2 3 5 9 14 15 | 2 5 6 9 11 15 | 3 4 7 11 14 15 | 4 5 7 8 12 13 | 5 9 10 13 14 15 |
| 1 2 4 12 13 15 | 1 3 6 7 9 11 | 1 5 7 8 11 13 | 2 3 5 10 13 14 | 2 5 6 12 14 15 | 3 4 7 12 13 14 | 4 5 7 9 12 14 | 6 7 8 9 10 14 |
| 1 2 5 6 7 9 | 1 3 6 7 10 12 | 1 5 7 9 11 15 | 2 3 5 11 12 15 | 2 5 7 8 9 15 | 3 4 8 9 12 13 | 4 5 7 9 13 15 | 6 7 8 9 11 12 |
| 1 2 5 6 8 12 | 1 3 6 8 10 15 | 1 5 7 10 13 14 | 2 3 6 7 8 12 | 2 5 7 8 10 14 | 3 4 8 10 11 13 | 4 5 7 10 11 12 | 6 7 8 9 11 13 |
| 1 2 5 6 11 12 | 1 3 6 8 11 12 | 1 5 7 11 12 14 | 2 3 6 7 9 10 | 2 5 7 8 11 15 | 3 4 8 10 12 15 | 4 5 7 10 11 14 | 6 7 8 9 11 15 |
| 1 2 5 6 13 14 | 1 3 6 8 11 13 | 1 5 8 9 14 15 | 2 3 6 7 11 14 | 2 5 7 8 12 15 | 3 4 8 11 12 14 | 4 5 7 10 11 15 | 6 7 8 10 13 15 |
| 1 2 5 7 9 13 | 1 3 6 9 13 14 | 1 5 8 10 11 12 | 2 3 6 8 10 13 | 2 5 7 9 10 11 | 3 4 8 13 14 15 | 4 5 8 9 10 14 | 6 7 9 10 11 14 |
| 1 2 5 7 10 12 | 1 3 6 10 11 14 | 1 5 8 10 12 14 | 2 3 6 8 11 15 | 2 5 7 9 12 14 | 3 4 9 11 12 13 | 4 5 8 9 11 12 | 6 7 11 12 13 14 |
| 1 2 5 7 14 15 | 1 3 6 11 12 13 | 1 5 8 12 13 15 | 2 3 6 9 11 12 | 2 5 7 11 12 13 | 3 4 9 12 14 15 | 4 5 8 10 13 15 | 6 8 10 11 13 15 |
| 1 2 5 8 9 10 | 1 3 7 8 10 14 | 1 5 9 10 11 14 | 2 3 6 9 13 14 | 2 5 8 9 11 12 | 3 4 10 11 12 14 | 4 5 8 11 13 14 | 6 8 10 12 13 15 |
| 1 2 5 8 11 14 | 1 3 7 8 11 14 | 1 5 9 10 12 13 | 2 3 6 10 11 13 | 2 5 8 9 13 14 | 3 5 6 7 8 15 | 4 5 8 12 14 15 | 6 8 10 13 14 15 |
| 1 2 5 9 12 15 | 1 3 7 8 12 14 | 1 5 9 12 13 14 | 2 3 6 10 14 15 | 2 5 8 10 12 13 | 3 5 6 7 9 12 | 4 5 9 10 11 13 | 6 9 10 11 12 13 |
| 1 2 5 10 11 13 | 1 3 7 8 14 15 | 1 5 10 11 13 15 | 2 3 7 8 9 14 | 2 5 8 10 14 15 | 3 5 6 7 11 14 | 4 5 9 11 14 15 | 6 9 10 12 14 15 |
| 1 2 5 10 13 15 | 1 3 7 9 10 12 | 1 5 11 12 14 15 | 2 3 7 8 10 15 | 2 5 9 11 13 14 | 3 5 6 8 10 14 | 4 5 9 12 13 15 | 6 9 11 13 14 15 |
| 1 2 5 11 13 15 | 1 3 7 9 10 15 | 1 6 7 8 12 15 | 2 3 7 8 12 13 | 2 5 10 11 12 14 | 3 5 6 9 11 13 | 4 5 10 12 13 14 | 7 8 9 12 13 14 |
| 1 2 6 7 8 10 | 1 3 7 11 13 15 | 1 6 7 8 13 14 | 2 3 7 9 11 13 | 2 5 10 11 13 15 | 3 5 6 9 12 15 | 4 6 7 8 11 14 | 7 8 10 11 12 14 |
| 1 2 6 7 11 13 | 1 3 7 12 13 15 | 1 6 7 9 10 13 | 2 3 7 10 11 15 | 2 5 12 13 14 15 | 3 5 6 10 11 12 | 4 6 7 8 12 13 | 7 8 10 12 14 15 |
| 1 2 6 7 12 14 | 1 3 8 9 12 13 | 1 6 7 9 12 14 | 2 3 7 12 14 15 | 2 6 7 8 9 11 | 3 5 6 10 13 15 | 4 6 7 9 10 12 | 7 8 11 12 14 15 |
| 1 2 6 8 9 13 | 1 3 8 9 12 15 | 1 6 7 10 11 15 | 2 3 8 9 10 11 | 2 6 7 8 14 15 | 3 5 6 13 14 15 | 4 6 7 10 11 13 | 7 9 10 11 13 15 |
| 1 2 6 8 11 15 | 1 3 8 10 13 15 | 1 6 7 11 14 15 | 2 3 8 9 12 15 | 2 6 7 9 12 13 | 3 5 7 8 9 10 | 4 6 7 10 14 15 | 7 10 11 12 14 15 |
| 1 2 6 9 10 11 | 1 3 9 10 11 13 | 1 6 8 9 10 12 | 2 3 8 10 11 12 | 2 6 7 9 12 15 | 3 5 7 8 10 13 | 4 6 7 12 14 15 | 8 9 10 11 12 13 |
| 1 2 6 9 11 14 | 1 3 9 11 14 15 | 1 6 8 9 14 15 | 2 3 8 10 11 14 | 2 6 7 9 13 15 | 3 5 7 8 11 12 | 4 6 7 13 14 15 | 8 9 12 13 14 15 |
| 1 2 6 10 12 13 | 1 3 10 11 12 15 | 1 6 8 10 11 14 | 2 3 8 13 14 15 | 2 6 7 10 11 12 | 3 5 7 9 11 15 | 4 6 8 9 10 13 | 8 10 11 12 14 15 |
| 1 2 6 10 14 15 | 1 3 12 13 14 15 | 1 6 8 11 12 13 | 2 3 9 10 12 14 | 2 6 7 10 13 14 | 3 5 7 10 14 15 | 4 6 8 9 12 15 | 9 10 11 12 13 14 |
| 1 2 7 8 9 14 | 1 4 5 6 7 13 | 1 6 9 10 13 15 | 2 3 9 11 13 15 | 2 6 7 12 13 15 | 3 5 7 12 13 14 | 4 6 8 9 13 14 | |
| 1 2 7 8 11 12 | 1 4 5 6 8 14 | 1 6 9 11 12 15 | 2 3 10 12 13 15 | 2 6 8 9 10 15 | 3 5 8 9 10 13 | 4 6 8 10 11 12 | |
| 1 2 7 8 13 15 | 1 4 5 6 9 12 | 1 6 10 11 12 14 | 2 3 11 12 13 14 | 2 6 8 9 12 14 | 3 5 8 9 12 14 | 4 6 8 11 14 15 | |
| 1 2 7 9 10 15 | 1 4 5 6 10 11 | 1 6 10 11 13 14 | 2 4 5 6 7 12 | 2 6 8 10 12 14 | 3 5 8 10 11 15 | 4 6 9 10 11 15 | |
| 1 2 7 9 11 12 | 1 4 5 6 11 14 | 1 6 12 13 14 15 | 2 4 5 6 8 10 | 2 6 8 11 12 14 | 3 5 8 11 14 15 | 4 6 9 11 12 14 | |
| 1 2 7 10 11 14 | 1 4 5 7 8 14 | 1 7 8 9 10 11 | 2 4 5 6 9 10 | 2 6 8 12 13 14 | 3 5 8 12 13 15 | 4 6 10 12 13 14 | |